\documentclass[11pt]{article}
\usepackage{amsmath,latexsym,amscd, amssymb, amsfonts}
\usepackage{graphics,amsthm}
\usepackage[dvips]{epsfig}
\newcommand{\RR}{\Bbb R}

\newcommand{\ZZ}{\Bbb Z}

\newcommand{\e}{\epsilon}
\newcommand{\para}{\kern0.2em{\backslash} \kern-0.7em {\backslash} \kern0.2em }
\newcommand{\parap}{          {\backslash} \kern-0.7em {\backslash} \kern0.2em }

\newcommand{\te}{\text{\textit{tub}}^{\epsilon}}

\newcommand{\CH}{{\mathcal H}}
\newcommand{\tr}{\times \RR}
\newcommand{\hCH}{\hat{{\mathcal H}}}
\newcommand{\hM}{\hat{M}}
\newcommand{\htheta}{\hat{\theta}}

\begin{document}

\newtheorem{thm}{Theorem}[section]
\newtheorem{prop}[thm]{Proposition}
\newtheorem{lemma}[thm]{Lemma}
\newtheorem{cor}[thm]{Corollary}
\newtheorem{claim}[thm]{Claim}
\newtheorem{fact}[thm]{Fact}

\newtheorem{defi}[thm]{Definition}
\newtheorem{pdef}[thm]{Proposition-Definition}
\newtheorem{ep}[thm]{Example}
\newtheorem{eps}[thm]{Examples}

\newtheorem{remark}[thm]{Remark}

\newcommand {\comment}[1]{{\marginpar{*}\scriptsize{\bf Comments:}\scriptsize{\ #1 \ }}}

\title {Averaging of Legendrian submanifolds of contact manifolds}
%Averaging of isotropic submanifolds and Weinstein's submanifold averaging}
\author{by \begin{Large}Marco Zambon\end{Large}\\
%Department of Mathematics, UC Berkeley, CA, USA\\
\texttt{zambon@math.unizh.ch}}
\date{\today}
\maketitle

\begin{abstract}
We give a procedure to ``average'' canonically $C^1$-close
Legendrian submanifolds of contact manifolds. As a corollary we
obtain that, whenever a compact group action leaves a Legendrian
submanifold almost invariant, there is an invariant Legendrian
submanifold nearby.
\end{abstract}

\section {Introduction}

Weinstein \cite{we} gave a construction which, given a
parameterized family $\{N_g\}$ of $C^1$-close submanifolds of a
Riemannian manifold $M$, produces an average submanifold $N$. The
construction depends only on the metric of $M$, therefore the
averaging procedure is equivariant with respect to isometries of
$M$.

 In \cite{za} we
adapted the above Riemannian averaging to the setting of
symplectic geometry: starting with a family of isotropic
submanifolds $\{N_g\}$ we construct canonically an isotropic
average submanifold. Applications include statements about almost
invariant isotropic submanifolds under group actions, images of
such submanifolds under moment maps, and equivariant
symplectomorphisms (see
 the introduction of \cite{za}).\\

In this short paper we consider the setting of contact
geometry\footnote{We would like to thank Xiang Tang for suggesting
this application to us and Alfonso Gracia-Saz, Tam\'as K\'alm\'an
and Alan Weinstein for useful comments.} and give a construction
to average Legendrian submanifolds.

Recall that a  \emph{contact manifold} is a manifold $M^{2n+1}$
together with a hyperplane distribution $\CH$ on $M$ such that
locally $\CH=\ker \theta$ for some locally defined 1-form $\theta$
satisfying
 $(d\theta)^n\wedge \theta \neq 0$.
A submanifold $N$ of $(M^{2n+1},\CH)$ is called Legendrian if it
is tangent to $\CH$ and it has maximal dimension among
submanifolds with this property, i.e. $\dim(N)=n$.

Consider the real line bundle $\CH^{\circ}\rightarrow M$
consisting of all covectors $\xi\in T^*M$ that annihilate the
distribution $\CH$. The distribution $\CH$ is co-orientable iff
this line bundle is trivial, and in this case there exists there
exists a \emph{contact one-form} representing $\CH$, i.e. a
(global) 1-form $\theta$ with kernel $\CH$ such that
$(d\theta)^n\wedge \theta$ is a volume form. The unique vector
field $E$ satisfying $\theta(E)=1$, $d\theta(E,\cdot)=0$ is called
Reeb vector field.

If the line bundle $\CH^{\circ}\rightarrow M$ is not trivial,
consider the manifold $\hat{M}$ of unit length (w.r.t. some
metric) elements of $\CH^{\circ}$. With the obvious projection
$\hat{M}$ is a (connected) double cover of $M$, and $(\hM,\hCH)$
is a co-orientable contact manifold, where $\hCH$ is the pullback
of the contact distribution $\CH$ to $\hM$.

\section{Results in the co-orientable case}\label{co}

In this section we state our results for a co-orientable contact
manifold $(M,\CH)$. In Section \ref{non-co} we will prove
analogous statements in the non co-orientable case by
 reducing it to the co-orientable one.

 First we need to
introduce some definitions from \cite{we}.

 If $M$ is a Riemannian manifold and $N$ a submanifold,
$(M,N)$ is called a \textit{gentle pair} if (i) the normal
injectivity radius of $N$ is at least 1; (ii) the sectional
curvatures of $M$ in the tubular neighborhood of radius one about
$N$ are bounded in absolute value by 1; (iii) the injectivity
radius of each point of the above neighborhood is at least 1.

The \textit{$C^0$-distance} between any two submanifolds $N$, $N'$
is denoted by $d_0(N,N')$ and is the supremum as $x'$ ranges over
$N'$ of the distance from $x'$ to the closest point in $N$.

The distance between two subspaces of the same dimension $F,F'$ of
a Euclidean vector space $E$, denoted by $d(F,F')$, is equal to
the $C^0$-distance
 between the unit spheres of
$F$ and $F'$ considered as Riemannian submanifolds of the unit
sphere of $E$, and is at most $\frac{\pi}{2}$. The
\textit{$C^1$-distance between two submanifolds} $N, N'$ of a
Riemannian manifold is defined whenever $N'$ lies in the tubular
neighborhood of $N$ and is the image under the exponential map of
a section of the normal bundle $\nu N$. It is denoted by $d_1(N,
N')$ and it is the supremum  as $x'$ ranges over $N'$
  of the length of the geodesic segment from $x'$ to
the nearest point $x$ in $N$ and the distance between $T_{x'}N'$
and the parallel translate of $T_xN$ along the above geodesic
segment.\\

Now consider a manifold with a contact form $(M,\theta)$, and
endow it with a \emph{compatible} Riemannian metric $g$ as
follows: for each fiber $\CH_p$ of the vector bundle $\CH=\ker
\theta \rightarrow M$, $(\CH_p,d\theta|_{\CH_p})$ is a symplectic
vector space, and we can choose a compatible positive inner
product $g$ (i.e. $d\theta(X,IY)=g(X,Y)$ determines an
endomorphism $I$ of $\CH_p$ satisfying $I^2=-Id$). We can do so in
smooth way (see \cite{ca}, Ch. 12). We extend $g$ to a Riemannian
metric on $M$ by imposing
that the Reeb vector field $E$ have unit length and be orthogonal to $\CH$.\\

We state our theorem for Legendrian submanifolds, even though it
 equally applies to submanifolds tangent to $\CH$ of lower
dimension. As a technical assumption we will require that the
$C^0$-norms of the covariant derivatives of $\theta$ and $d\theta$
with respect to the Levi-Civita connection be bounded by $1$ (but
see Remark \ref{bound}).

\begin{thm}\label{leg}
Let $(M,\theta)$ be a manifold with a contact form, endowed with a
Riemannian metric $g$ as above so that $|\nabla \theta|, |\nabla
d\theta|<1$. Let $\{N_g\}$ be a family of Legendrian submanifolds
of $M$ parameterized in a measurable way by elements of a
probability space $G$, such that all the pairs $(M,N_g)$ are
gentle. If $d_1(N_g,N_h)<\epsilon<\frac{1}{70000}$ for all $g$ and
$h$ in $G$, there is a well defined \emph{\textbf{Legendrian
average}} submanifold $L$ with $d_0(N_g,L)<1000 \e$ for all $g$ in
$G$. This construction is equivariant with respect to isometric
contactomorphisms of $(M,\theta)$ and measure preserving
automorphisms of $G$.
\end{thm}

\begin{remark} \label{bound}
The theorem holds even if the bound on $|\nabla \theta|$ and
$|\nabla d\theta|$ is larger than $1$, but in that case the bound
on $\epsilon$ has to be chosen smaller.
\end{remark}

A simple consequence, which we want to state in terms of contact
manifolds, is the following:
\begin{thm}\label{gr} Let $(M,\CH)$ be a co-orientable
contact manifold, let $G$ be a compact Lie group acting on $M$
preserving $\CH$ and its co-orientation, and let ${N_0}$ be a
Legendrian submanifold. Endow $(M,\CH)$ with a contact form
$\theta$ and a Riemannian metric $g$ as above, both invariant
under the  $G$ action. Suppose that $|\nabla \theta|, |\nabla
d\theta|<1$. Then if $(M,N_0)$ is a gentle pair and
$d_1({N_0},g{N_0})<\e<\frac{1}{70000}$ for all $g \in G$, there
exists a $G$-invariant Legendrian submanifold $L$ of $(M, \CH)$
with $d_0({N_0},L)< 1000\e$.
\end{thm}

Indeed, we just need to endow $G$ with its bi-invariant probability measure and apply Theorem \ref{leg} to the family $\{g N_0\}$:
 their average will be $G$-invariant by the equivariance properties of Theorem \ref{leg}.

\begin{remark}\label{rem} We can always find $\theta$ and $g$ which are
$G$-invariant:
 by averaging over $G$ we can obtain a $G$-invariant one-form $\theta$ representing $\CH$,
and  using some $G$-invariant  metric on the vector bundle $\CH \rightarrow M$ as a tool
 we can construct a ``compatible''  metric $g$ which is $G$-invariant  (see Ch. 12 in \cite{ca}). However in general we can not give any a priori bounds on
the covariant derivatives of $\theta$ and $d\theta$.
\end{remark}
\begin{remark}
The construction of Theorem \ref{leg} is actually equivariant also
w.r.t. isometric diffeomorphisms of $M$ that preserve $\theta$ up
to a sign (see the proof of Theorem \ref{legnon}). Using this and
Remark \ref{invnon} below one sees that Theorem \ref{gr} holds
even without assuming that the $G$ action preserve the
co-orientation of $\CH$.
\end{remark}

\section{Idea of the proof of Theorem \ref{leg}}

An approach to prove Theorem \ref{leg} is to use the idea that
worked in the symplectic setting (see Section 1.3 in \cite{za}).
This would be carried out as follows: first construct the
Weinstein average $N$ of the Legendrian submanifolds $\{N_g\}$. For each $g\in G$, using the
Riemannian metric construct a diffeomorphism $\varphi_g$ from a
neighborhood of $N_g$ to a neighborhood of $N$ and denote by
$\theta_g$ the pullback of $\theta$ to the neighborhood of $N$.
Since the submanifolds $N_g$ are close to $N$, each form in the convex
linear combination $\theta_t=\theta +t(\int \theta_g -\theta)$ is
a contact form, say with contact distribution $\CH_t$. Therefore
we can apply the contact version of Moser's Theorem (see
\cite{ca}, Ch. 10). It states that if the vector field $v_t$ is the inverse image of
$-(\int \theta_g -\theta)|_{\CH_t}$  by (the isomorphism induced
by) $d\theta_t|_{\CH_t}$, then the time-one flow $\rho_1$ of
$\{v_t\}$ satisfies $(\rho_1)_*\CH=\CH_1$. Therefore, since $N$ is
tangent to $\CH_1$, its pre-image $L$ under $\rho_1$ is tangent to
$\CH$, i.e. it is a Legendrian submanifold of $(M,\theta)$.

This construction can indeed be carried out and satisfies the
invariance properties stated in Theorem \ref{leg} since all steps are
canonical. However it delivers a numerically quite unsatisfactory
 estimate for $d_0(N_g,L)$, therefore we choose \emph{not} to use
this  approach but rather a different one, which we outline now.\footnote{The estimates needed for our first approach are completely analogous
 to those needed for the second approach, i.e. those of \cite{za}. A difference though is that in the first approach we
make use of a  bound on the norm of the $C^0$-small one-form
$\int_g \theta_g -\theta$, whereas in the second approach we will
need the norm of a primitive of a certain $C^0$-small two-form.
While passing to the primitive we will improve the $C^0$-norm by
``one order of magnitude'' (see equation ($\bigstar$) in Section
7.2 of \cite{za}), and this is responsible for the better estimates
obtained using the second approach.}\\

Recall that the \emph{symplectization} of a manifold with contact form
$(M,\theta)$ is the symplectic manifold $(M \times \RR, d(e^s
\theta))$, where $s$ denotes the coordinate on the $\RR$ factor.
Here and in the following we abuse notation by writing $\theta$ in
place of $\pi^* \theta$, where $\pi: M \tr \rightarrow M$.
   The  main observation is the
following lemma, whose proof consists of a short computation and is
omitted.
\begin{lemma}
$L$ is a Legendrian submanifold of $M$ if and only if $L\times
\RR$ is a Lagrangian submanifold of $M \times \RR$.
\end{lemma}

We make use of the averaging procedure for Lagrangian submanifolds
(\cite{za}), which reads:

\begin{thm} \label{iso} Let $(\tilde{M},g,\omega)$ be an almost-K\"ahler
manifold\footnote{This means that the symplectic form $\omega$
and the metric $g$ are compatible.}
 satisfying $|\nabla \omega|<1$ and $\{\tilde{N}_g\}$ a family
of Lagrangian submanifolds of $\tilde{M}$ parameterized in a
measurable way by elements of a probability space $G$, such that
all the pairs $(\tilde{M},\tilde{N}_g)$ are gentle. If
$d_1(\tilde{N}_g,\tilde{N}_h)<\epsilon<\frac{1}{70000}$ for all
$g$ and $h$ in $G$, there is a well defined
\emph{\textbf{Lagrangian average}} submanifold $\tilde{L}$ with
$d_0(\tilde{N}_g,\tilde{L})<1000 \e$ for all $g$ in $G$. This
construction is equivariant with respect to isometric
symplectomorphisms of $\tilde{M}$ and measure preserving
automorphisms of $G$.\end{thm}

Now the strategy for our (second) approach is straightforward: given the
family $\{N_g\}$ of Legendrian submanifolds of $(M,\theta)$ we
 consider the Lagrangian family $\{N_g \times \RR\}$ in the
symplectization of $M$, apply the Lagrangian averaging theorem,
and \emph{if} the average is a product $L \times \RR$ then $L$
will be our Legendrian average. The invariance properties stated in Theorem  \ref{leg} are satisfied because this construction is canonical after choosing the contact form and metric on $M$.

\section{The proof of Theorem \ref{leg}}

We endow the symplectization $M \tr$ with the product metric
obtained from $(\RR, ds \otimes ds)$ and $(M,g)$, and by abuse of
notation denote this metric by $g$. Unfortunately $(M \tr, d(e^s
\theta), g)$ does \emph{not} satisfy the assumptions of Theorem
\ref{iso}: indeed it is not an almost-K\"ahler manifold  (however
$g$ is compatible with the non-degenerate 2-form $e^{-s}d(e^s
\theta)= ds \wedge \theta + d\theta$). Furthermore the condition
on the boundedness of $\nabla d(e^s \theta)$ is also violated.
Therefore we can not just apply Theorem \ref{iso} but we have to
follow the construction in the theorem and check that it applies
to $(M \tr, d(e^s \theta), g)$.

The remaining assumptions of Theorem \ref{iso} concerning the
gentleness of the pairs $(M \tr, N_g \tr)$ and the distances
$d_1(N_g \tr, N_h \tr)$ are satisfied, since the metric $g$ on $M
\tr$ is a product metric.\\

In the remainder of the proof we will follow the construction of
Theorem \ref{iso}. We do so for two reasons: firstly, in order to
make sure that the Lagrangian average of the $\{N_g \times \RR\}$
is of the form $L \times \RR$ for some submanifold $L$ of $M$, and
secondly
  to check that the construction applies to $(M \tr,
\omega:=d(e^s \theta), g)$ even though the assumptions of Theorem
\ref{iso} are not satisfied.

We refer the reader to Section 1.3 of \cite{za} for the outline
of the construction of Theorem \ref{iso} and adopt the notation
used there.
 We divide the construction
into 5 steps.\\

\textbf{Step 1:} \emph{Construct the Weinstein average of the
family $\{N_g \tr\}$}.

\noindent Step 1 applies to our manifold $M \tr$ because
 Step 1
involves only the Riemannian structure of $M \tr $, which
satisfies the assumptions of Weinstein's averaging theorem (see
\cite{za} and \cite{we}).
 Notice that
the Weinstein average is of the form $N \tr$ because the group
$\RR$ acts isometrically on $M \tr$ by translation of the second
factor and Weinstein's averaging
procedure is equivariant w.r.t. isometries.\\

 \textbf{Step 2:} \emph{The restriction of
$$\varphi_g: \text{ \emph{Neighborhood of }}N_g \tr \rightarrow
\text{\emph{Neighborhood of }}N \tr $$ to $\exp_{N_g \tr} (\nu
(N_g \tr))_{0.05}$ is a diffeomorphism onto. }

\noindent Step 2 applies to $M \tr$ because it involves only its
Riemannian structure. Notice that $\varphi_g$ preserves $M \times
\{s\}$ for each $s \in \RR$ and that the $\varphi_g|_{M \times
\{s\}}$ coincide for all $s$ (under the obvious identifications $M
\times \{s\} \equiv M \times
\{s'\}$ ), since the metric on $M \tr$ is the product metric.\\

 \textbf{Step 3:} \emph{ On $\te$
the family $\omega_t:=\omega+t(\int_g \omega_g -\omega)$ consists
of symplectic forms lying in the same cohomology class, where
$\omega_g:= (\varphi_g^{-1})^*\omega$ and $\te$ is neighborhood of $N$.
 }

\noindent
% Because of the remark in Step 2 we have
%$\omega_g=e^s\cdot (\varphi_g^{-1})^*\bar{\omega}= e^2 \cdot
%\bar{\omega_g}$ where
Define $\bar{\omega}:=ds \wedge \theta + d\theta$, a
non-degenerate 2 form compatible with the metric $g$ on $M \tr$.
Below we will show that $|\nabla_X \bar{\omega}| \le 2$ for
``horizontal'' unit vectors $X$. Since $\omega=e^s\bar{\omega}$ using
this we see that the statements of Lemma 7.2 and Corollary 7.1 in
\cite{za} hold if one multiplies the right hand sides there by
$e^s$ and multiplies by 2 the term ``$2L+100 \epsilon$''\footnote{This factor of 2 does not affect the numerical estimates that follow.}.
This shows that the $\omega_t$ are non-degenerate on $\te$. One
sees that the $\omega_t$ lie in the zero cohomology class exactly
as in Section 1.3 of \cite{za}.

Now we derive the bound on $|\nabla_X \bar{\omega}|$, where $X \in
T(M \times \{s\})$ is a unit vector:
\begin{eqnarray*}
|\nabla_X (ds \wedge  \pi^*\theta + d \pi^*\theta)| &=& |\nabla_X
ds\wedge \pi^* \theta + ds \wedge \nabla_X \pi^* \theta + \nabla_X
\pi^* d\theta|\\
&=& |ds|\cdot |\nabla_X\pi^* \theta|+|\nabla_X \pi^* d\theta| \\
&\le& 1\cdot  1+1.
\end{eqnarray*}

 \textbf{Step 4:} \emph{Construct canonically
 a primitive $\alpha$ of $\int_g \omega_g -\omega$. }

\noindent A primitive of $\omega_g -\omega$ is given in a
canonical fashion by $e^s( (\varphi_g^{-1})^* \theta -\theta)$,
however we do \emph{not} want to use this primitive since it would
deliver bad numerical estimates, as it happens in our first
approach. Instead we construct one using the (more involved)
procedure of Section 7.2 in \cite{za}.

 Because of the remark in Step 2 we can write
$\omega_g=e^s \bar{\omega}_g$ where
$\bar{\omega}_g:=(\varphi_g^{-1})^*\bar{\omega}$, therefore the
2-form $\omega_g -\omega$
 can be written as $e^s(\bar{\omega}_g -\bar{\omega})$. Notice that $\bar{\omega}$ is ``constant
in $\RR$-direction'' (i.e. the Lie derivative
$L_{\frac{\partial}{\partial s}}
 \bar{\omega}=0$), and by the
remarks in Step 2 the same holds of $\bar{\omega}_g$.
 The construction in \cite{za} commutes with multiplication of forms by
functions of $s$  and furthermore preserves the condition of being
``constant in $\RR$ direction'' (this follows from the explicit
formula for $\alpha^g$ in Section 7.2 of \cite{za} and from the
fact that the vector bundles used there are pullbacks of vector
bundles over subsets of $M$ via  $\pi:M \times \RR \rightarrow
M$). Therefore that construction applied to $\int \omega_g
-\omega$ delivers a primitive $\alpha= e^s \bar{\alpha}$ where
$\bar{\alpha}$ is ``constant in $\RR$ direction''. Furthermore
$|\bar{\alpha}|$ is estimated as in Proposition 7.2 of \cite{za},
because $\bar{\omega}$ is compatible
with the metric.\\

 \textbf{Step 5:} \emph{Obtain the Lagrangian average
by following backwards the Moser vector field
 $v_{t}:=-\omega_t^{-1}(\alpha)$ starting from $N \tr$.}

\noindent The key observation is that
$$v_t=
-\omega_t^{-1}(\alpha)=-\bar{\omega}_t^{-1}(\bar \alpha),$$
therefore the vector field $v_t$ is independent of $s$ (i.e.
$L_{\frac{\partial}{\partial s}}v_t=0$ for $t\in [0,1]$). This
implies that the Lagrangian average is of the form $L \times \RR$.
Further, this implies that $|v_t|$ is bounded as in Proposition
7.3 of \cite{za}, so that the proof of the theorem can be
concluded verbatim as in \cite{za}.

\section{Results in the non co-orientable case}\label{non-co}

Until the end of this note  $(M,\CH)$ will be a \emph{non}
co-orientable contact manifold. Even though $M$ can not be endowed
with a (global) contact form representing $\CH$, there is a
``1-form which is well-defined up to a sign'' representing $\CH$,
as follows. Recall from the introduction that there is a
co-orientable double cover $\pi:(\hM,\hCH)\rightarrow (M,\CH)$.
There is an obvious $\ZZ_2$ action on $\hM$; denote by
$i:\hM\rightarrow \hM$ the action of the non-trivial element of
$\ZZ_2$. Since we assume that $\CH$ is not co-orientable, $i$
 inverts
the co-orientation of $\hCH$, hence if we choose any contact form
$\theta$ on $\hM$ with kernel $\hCH$, the average
$\htheta:=\frac{1}{2}(\theta-i^*\theta)$ will be again a contact
form. This contact form satisfies $i^*\htheta=-\htheta$, so it
does not descend to $M$. However it gives rise to a subset of
$T^*M$, namely
%$\{ \pm
%\theta_p|p\in \hat{M}\}$ of $T^*\hM$ which descends via
%$\pi:\hM\rightarrow M$ to a subset
$\Theta= \{\xi\in T_{p}^*M: \pi^*(\xi)=\htheta_{\hat{p}} \text{
for a } \hat{p}\in \pi^{-1}(p)\}$.
%Notice that the restriction of the
%canonical projection $T^*M\rightarrow M$ to this subset is a
%(connected) double covering and that $\Theta$ is invariant under
%multiplication my $-1$ in $T^*M$. In other words,
Notice that on a small open subset of $M$,
%under $\pi:\Theta\rightarrowM$ corresponds
$\Theta$ corresponds to two 1-forms on that open set, which are
one the negative of the other, and which pull back via $\pi$ to
$\pm \htheta$.
 A metric
compatible with $\Theta$ is one that on small open sets is
compatible with one (or equivalently both) of these two 1-forms;
such a metric always exists.
%it arises by taking a $\ZZ_2$-invariant metric compatible to
%$\theta$ on $\hM$ and pushing it down to $M$.
%(Such a metric can always be constructed by choosing a $\ZZ_2$
%invariant metric $T\hM$ and making it compatible to $\htheta$ as
%in Remark \ref{rem} above).
In order to state the next theorem in terms on $M$ alone, we abuse
notation by defining $|\nabla \Theta|$ to be supremum (over all
small open sets in $M$) of the $C^0$-norm of the covariant
derivative of any of the two local 1-forms that locally correspond
to $\Theta$, and similarly for $|\nabla d\Theta|$.

Using the above notation we have:
\begin{thm}\label{legnon}
Let $(M,\CH)$ be a non co-orientable contact manifold. Choose
$\Theta\subset T^*M$ as above, and endow $M$ with a Riemannian
metric $g$ compatible with $\Theta$ so that $|\nabla \Theta|,
|\nabla d\Theta|<1$ on $\hM$. Let $\{N_g\}$ be a family of
Legendrian submanifolds of $(M,\CH)$ parameterized in a measurable
way by elements of a probability space $G$, such that all the
pairs $(M,N_g)$ are gentle. If
$d_1(N_g,N_h)<\epsilon<\frac{1}{70000}$ for all $g$ and $h$ in
$G$, there is a well defined \emph{\textbf{Legendrian average}}
submanifold $L$ with $d_0(N_g,L)<1000 \e$ for all $g$ in $G$. This
construction is equivariant with respect to isometries of $M$ that
preserve $\Theta$ and measure preserving automorphisms of $G$.
\end{thm}

\begin{proof}
The key point is that applying the construction of Theorem
\ref{leg} to contact manifolds $(\bar{M},\bar{\theta})$ and
$(\bar{M},-\bar{\theta})$ (with the same choice of compatible
metric) leads to the same averaging. Indeed, Steps 1 and 2 do not
involve the contact form, whereas the forms $\omega_t$ in Steps 3
differ by a sign. The same holds for the form  $\alpha$ in Step 4.
In the construction of the Moser vector field
$v_t=-\omega_t^{-1}(\alpha)$ these signs cancel out, and we obtain
the same Moser vector field on $\bar{M}$ independently of whether
we use $\bar{\theta}$ or $-\bar{\theta}$. Therefore we can improve
Theorem \ref{leg} by saying that the construction
 is equivariant with respect to isometries of $\bar{M}$ that preserve
$\bar{\theta}$ up to a sign.

Now the proof follows easily (again we refer to above for the
notation): the $\ZZ_2$ action on $\hM$ is by isometries (w.r.t.
the pullback metric)  and preserves $\htheta$ up to a sign. We can
lift the $\{N_g\}$'s to a family $\{\hat{N_g}\}$ of Legendrian
submanifolds of $(\hM,\htheta)$. Since the family $\{\hat{N_g}\}$
 is $\ZZ_2$ invariant,
by the above equivariancy
  their average is also $\ZZ_2$ invariant,
 hence it descends to a Legendrian submanifold $L$ of $M$.
 This averaging procedure is equivariant w.r.t. isometries $\phi$ of $M$
 that preserve (via the cotangent lift) $\Theta\subset T^*M$,
 since it depends only on $\Theta$ and the metric.
 %Indeed $\phi$ (again via the cotangent lift) induces a
 %diffeomorphism $\hat{\phi}$ of $\hM$, which preserves  $\htheta$, possibly
 %up to a  sign,  and of course the metric.
%Hence averaging the $\{\hat{N_g}\}$'s commutes with $\hat{\phi}$,
%and by projecting everything to $M$ we see that averaging the
%$\{{N_g}\}$'s commutes with ${\phi}$.
\end{proof}

The analog of Theorem \ref{gr} follows from the equivariance
properties of Theorem \ref{legnon}:
\begin{thm}\label{grnon} Let $(M,\CH)$ be a non co-orientable
contact manifold, let $G$ be a compact Lie group acting on $M$
preserving $\CH$, and let ${N_0}$ be a Legendrian submanifold.
Endow $(M,\CH)$ with $\Theta\subset T^*M$ as above and with a
Riemannian metric $g$ compatible with
 $\Theta$, both invariant under the  $G$ action.
Suppose that $|\nabla \Theta|, |\nabla d\Theta|<1$.
  Then if $(M,N_0)$ is a
gentle pair and $d_1({N_0},g{N_0})<\e<\frac{1}{70000}$ for all $g
\in G$, there exists a $G$-invariant Legendrian submanifold $L$ of
$(M, \CH)$ with $d_0({N_0},L)< 1000\e$.
\end{thm}

\begin{remark}\label{invnon}
We can always find $\Theta$ and $g$ as above. As seen at the
beginning of this section we can find a contact form $\htheta$ on
$\hM$ satisfying $i^*\htheta=-\htheta$. Define
$$\htheta_G:=\int_G \sigma(g)g^*\htheta,$$
where we are considering the lifted action of $G$ on $\hM$ and
$\sigma(g)$ is either $1$ or $-1$ depending on whether the action
of $g$ (and hence of the whole connected component of $G$ in which
$g$ lies) preserves or inverts the co-orientation of $\htheta$.
Being the average of contact forms defining the same
co-orientation of $\hCH$, $\htheta_G$ is a contact form, and
$h^*\htheta_G=\sigma(h) \htheta_G$.
%(with the plus or minus sign
%depending on whether the action of $h\in G$ on $\hM$ preserves the
%co-orientation of $\htheta$ or not).
 Further
$i^*\htheta_G=-\htheta_G$, so the $G$ action preserves the subset
$\Theta$ of $T^*M$ induced by $\htheta_G$. Modifying slightly
Remark \ref{rem} we can construct a metric on $\hM$ which is
invariant under the $G\times \ZZ_2$ action and compatible with
$\htheta_G$. This metric descends to a $G$ invariant metric on $M$
which is compatible with $\Theta$.
\end{remark}


\begin{thebibliography}{99}

%\bibitem [BK]{bk} Buser, P., Karcher, H.: Gromov's almost flat manifolds, Asterisque \textbf{81},
%Societe Mathematique de France, Paris, 1981\\
\bibitem [Ca]{ca} Cannas da Silva, A.: Lectures
on Symplectic Geometry. Springer Lecture Notes, Vol.
1764. Berlin, Heidelberg, New York: Springer 2000\\
%\bibitem [Jo]{jo} Jost, J.: Riemannian Geometry
%and Geometric Analysis, second edition. Berlin, Heidelberg, New York: Springer 1998\\
\bibitem [We]{we} Weinstein, A.: Almost
invariant submanifolds for compact group actions. J. Eur. Math. Soc. \textbf{2}, 53-86 (1999)\\
\bibitem [Za] {za} Zambon, M.: Submanifold averaging in riemannian and symplectic geometry, arXiv.DG/0208203\\

\end{thebibliography}
\end{document}